\input amstex
\input pstricks
\documentstyle{amsppt}
\magnification 1200
\newcount\defheadno
\newcount\headno
\newcount\subheadno
\newcount\formno
\newcount\cno
\newcount\figno
\global\formno=0
\global\headno=0

\def\nextsubheadno{\global\advance\subheadno by 1 \the\headno.\the\subheadno}
\def\nextheadno{\global\advance\headno by1 \the\headno \global\subheadno=0}
\def\nextno{\global\advance\cno by1 \the\cno }
\def\nextfigno{\global\advance\figno by1 \the\figno }
\def\nextformno{\global\advance\formno by1 \the\formno }

\def\blankpage{\nopagenumbers\eject\line{}\vfil\eject
\global\advance\pageno by -1\footline={\hss\tenrm\folio\hss}}
\def\eqnum{\tag{\nextformno}}
\def\eqlabel#1{\edef #1{\the\formno}}
\def\eqref#1{$(#1)$}

\def\figlabel#1{\edef #1{\the\figno}}

\def\headnum{\nextheadno. }
\def\headlabel#1{\edef #1{\the\headno}}
\def\headref#1{#1}
\def\subheadnum{\nextsubheadno. }
\def\subheadlabel#1{\edef #1{\the\subheadno}}

\def\procnum{\nextno}
\def\proclabel#1{\edef #1{\the\cno}}
\def\procref#1{$#1$}
\def\N{\Bbb N}

\def\Z{\Bbb Z}
\def\R{\Bbb R}

\def\E{\Bbb E}
\def\J{\Bbb J}

\def\B{\Bbb B}
\def\Phii{\Phi_\infty}

\def\tend{\rightarrow}

\def\pg{p_g}
\def\pgd{\widehat{p}_g}

\def\diam{\text{diam}}
\def\und#1{\underline{#1}}
\def\pr{\text{pr}}
\def\bl#1{\bold{#1}}
\def\brond{\smash{B \raise 7.9pt\hbox{$\kern -9.4pt{\scriptstyle \circ}$}}}
\let\Dsty\displaystyle

\newcount\refno
\global\refno=0
\def\nextrefno{\global\advance\refno by 1 }
\nextrefno\edef\AlexI{\the\refno}
\nextrefno\edef\AlexII{\the\refno}
\nextrefno\edef\AlexIII{\the\refno}
\nextrefno\edef\AlexIV{\the\refno}
\nextrefno\edef\AlexV{\the\refno}
\nextrefno\edef\ACC{\the\refno}
\nextrefno\edef\AntPisz{\the\refno}
\nextrefno\edef\Barlow{\the\refno}
\nextrefno\edef\Bodineau{\the\refno}
\nextrefno\edef\Cerf{\the\refno}
\nextrefno\edef\CerfPiszI{\the\refno}
\nextrefno\edef\CerfPiszII{\the\refno}
\nextrefno\edef\CCS{\the\refno}
\nextrefno\edef\DeuschPisz{\the\refno}
\nextrefno\edef\DH{\the\refno}
\nextrefno\edef\DKS{\the\refno}
\nextrefno\edef\Grimo{\the\refno}
\nextrefno\edef\GrimoII{\the\refno}
\nextrefno\edef\Grim{\the\refno}
\nextrefno\edef\GrimII{\the\refno}
\nextrefno\edef\GrimMarst{\the\refno}
\nextrefno\edef\GrimPiza{\the\refno}
\nextrefno\edef\Hryniv{\the\refno}
\nextrefno\edef\IoffeI{\the\refno}
\nextrefno\edef\IoffeII{\the\refno}
\nextrefno\edef\IoffSchon{\the\refno}
\nextrefno\edef\LSS{\the\refno}
\nextrefno\edef\MathieuRemy{\the\refno}
\nextrefno\edef\PenPisz{\the\refno}
\nextrefno\edef\Pfis{\the\refno}
\nextrefno\edef\PfisVelen{\the\refno}
\nextrefno\edef\Pisz{\the\refno}
\nextrefno\edef\SSI{\the\refno}
\nextrefno\edef\SSII{\the\refno}
\nextrefno\edef\SSIII{\the\refno}
\nextrefno\edef\SSIV{\the\refno}
\topmatter
\NoRunningHeads
\title
Surface order Large deviations for 2D FK-percolation and Potts  models
\endtitle
\author
Olivier Couronn\'e, Reda J\"urg Messikh
\endauthor
\affil
Universit\'e Paris-Sud
\endaffil
\address
   Universit\'e Paris-Sud. Laboratoire de math\'ematiques,
       b\^at. 425, 91405 Orsay Cedex, France.
 \endaddress
\email
couronne\@cristal.math.u-psud.fr,
reda-jurg.messikh\@math.u-psud.fr
\endemail
\thanks
We would like to thank R. Cerf for suggesting the problem and for many helpful discussions.
\endthanks
\date
10 March 2003
\enddate
\abstract
By adapting the renormalization techniques of Pisztora, \cite{\Pisz},
we establish surface order large deviations estimates for FK-percolation on $\Z^2$ with
parameter $q\geq 1$ and for the corresponding Potts models.
Our results are valid up to the exponential decay threshold of dual connectivities which is widely believed to agree with the critical point.
\endabstract
\endtopmatter

\document
\subhead{\headnum Introduction and statement of results}\endsubhead In this paper we derive
surface order large deviations for Bernoulli percolation,
FK-percolation with parameter $q>1$ and for the corresponding
Potts models on the planar lattice $\Z^2$.
\par In dimension two, surface order large deviations behaviour and the Wulff
construction has been established for the Ising model
\cite{\DH, \DKS, \Hryniv, \IoffeI, \IoffeII, \IoffSchon, \Pfis, \PfisVelen, \SSI, \SSII, \SSIII, \SSIV},
for independent percolation \cite{\AlexIII, \ACC} and for the random cluster model
\cite{\AlexIV}. These works include also more precise results than large deviations for
the Wulff shape. They are obtained by using the skeleton coarse
graining technique to study dual contours which represent the interface. In higher dimensions other
methods had to be used to achieve the Wulff construction,
\cite{\Bodineau, \Cerf, \CerfPiszI, \CerfPiszII}, where one of the main tools that have been used was the
blocks coarse graining of Pisztora \cite{\Pisz}. This renormalization
technique led to surface order large deviations estimates for FK-percolation and for the
corresponding Potts models simultaneously. The results of
\cite{\Pisz}, and thus the Wulff construction in higher dimensions, are valid up to the
limit of the slab percolation thresholds. In the case of independent percolation, this
threshold has been proved to agree with the critical point \cite{\GrimMarst} and is
believed to be so for all the FK-percolation models with parameter $q\geq 1$ in dimension
greater than two.
\par Our aim is to adapt Pisztora's techniques to the two-dimensional
lattice thereby opening the way to an other proof for the Wulff
construction in dimension two. It is also worth noting that Pisztora's
renormalization technique forms a building block that has been used to answer
various other questions related to percolation \cite{\AntPisz,
\Barlow, \MathieuRemy, \PenPisz}. Thus, we expect that adapting \cite{\Pisz} to the
two-dimensional case will permit the use of this building block for
other problems on the planar lattice.  The main point in
our task is to get rid of the percolation in slabs which is specific
to the higher dimensional case. For this we produce estimates analogue to those of theorem 3.1 in
\cite{\Pisz} relying on the hypothesis that the dual connectivities decay exponentially.
This hypothesis is very natural in $\Z^2$, because it is possible to translate events from the
supercritical regime to the subcritical regime by planar duality. For Bernoulli
percolation, the exponential decay of the connectivities is known to hold in all the
subcritical regime, see \cite{\Grimo} and the references therein. For the random cluster
model on $\Z^2$ with $q=2$ the exponential decay follows from the exponential decay of
the correlation function in the Ising model \cite{\CCS}, and a proof has also been given
when $q$ is sufficiently large, see \cite{\Grim} and the references therein. Even if not
proved, the exponential decay of the connectivities is widely believed to
hold up to the critical point of all the FK-percolation models with
$q\geq 1$. In addition to that, we use a property which is specific to
the two dimensional case, namely the weak mixing property. This
property has been proved to hold for all the random cluster models
with $q\geq 1$ in the regime where the connectivities decay
exponentially \cite{\AlexI}. We need this property in order to use the
exponential decay in finite boxes \cite{\AlexII}.
\subsubhead \subheadnum {Statement of results} \endsubsubhead
Our results concern asymptotics of FK-measures on finite boxes $B(n)=(-n/2,
n/2]^2\cap\Z^2$, where $n$ is a positive integer. We will denote by $\Cal R(p,q,B(n))$
the set of these FK-measures defined on $B(n)$ with parameters $(p,q)$ and where we have
identified some vertices of the boundary. For $q\geq 1$ and $0<p\not= p_c(q)<1$, it is
known \cite{\GrimII} that there is a unique infinite volume Gibbs measure that we will
note $\Phii^{p,q}$. It is also known that $\Phii^{p,q}$ is translation invariant and
ergodic. In the uniqueness region, we will denote by $\theta=\theta(p,q)$ the density of
the infinite cluster. As the exponential-decay plays a crucial rule in our analysis, we
will introduce the following threshold \footnote{The notation $p_g$ comes from
\cite{\Grim}.}
 $$
\pg =\sup\{p:\exists c>0,\,\forall\; x\  \forall\; y\in \Z^2,\
\Phi^{p,q}_{\infty}[ x\leftrightarrow y]\leq \exp(-c|x-y|)\},
\eqnum
 $$\eqlabel{\expthreshold}\noindent
where $|x-y|$ is the $\Cal L^1$ norm and $\{x\leftrightarrow y\}$ is the event that there
exists an open path joining the vertex $x$ to the vertex $y$.
\par\noindent By the results of \cite{\GrimPiza}, it is known that exponential decay holds
as soon as the connectivities decay at a sufficient polynomial rate. We thus could replace
\eqref{\expthreshold} by
 $$
\pg =\sup\{p:\exists c>0,\,\forall\; x\  \forall\; y\in \Z^2,\
\Phi^{p,q}_{\infty}[ x\leftrightarrow y]\leq c/|x-y|)\}.
 $$
We introduce the point dual to $\pg$:
$$
\pgd=\frac{q(1-\pg)}{\pg+q(1-\pg)}\geq p_c(q),
$$
which is conjectured to agree with the critical point $p_c(q)$.
\par Our result states that up to large deviations of surface order,
there exists a unique biggest cluster in the box $B(n)$ with the same
density than the infinite cluster, and that the set of clusters of
intermediate size have a negligible volume. To be more precise, we
say that a cluster in $B(n)$ is {\it crossing} if it intersects all the faces of $B(n)$.
For $l\in\N$, we say that a cluster is $l$-{\it intermediate} if it is not of maximal
volume and its diameter does exceed $l$. We denote by $\J_l$ the set of $l$-intermediate
clusters. Let us set the event
 $$\eqalign{
 K(n,\varepsilon,l)=&\Big\{\exists!\hbox{ open cluster }C_m \hbox{ in }B(n)\hbox{ of
 maximal volume, }
C_m \hbox{ is crossing}\Big\}\cr
&\cap\Big\{
 n^{-2}|C_m|\in(\theta-\varepsilon,\theta+\varepsilon)\Big\}\cap
\Big\{ n^{-2}\sum_{C\in\J_l}|C|<\varepsilon\Big\},
}$$

\proclaim{Theorem \procnum}
Let $q\geq 1, 1>p>\pgd$ and $\varepsilon\in(0,\theta/2)$ be fixed. Then there exists a
constant $L$ such that
$$
\eqalign{ -\infty<&
\liminf_{n\rightarrow\infty}\frac{1}{n}\log\inf_{\Phi\in\Cal R(p,q,B(n))}
\Phi[K(n,\varepsilon,L)^c] \cr
\leq &
\limsup_{n\rightarrow\infty}\frac{1}{n}\log\inf_{\Phi\in\Cal R(p,q,B(n))}
\Phi[K(n,\varepsilon,L)^c]<0.}
$$\endproclaim\proclabel{\prin}
\noindent This result, via the FK-representation, can be used as in \cite{\Pisz}
to deduce large deviations estimates for the magnetization of the
Potts model. We will omit this as it would be an exact repetition of
theorem $1.1$ and theorem $5.4$ in \cite{\Pisz}.
\subsubhead{\subheadnum Organization of the paper}\endsubsubhead
In the following section we introduce notation and give a summary of
the FK model and of the duality in the plane. In section~4, we study
connectivity properties of FK percolation in a large box $B(n)$ and
establish estimates that will be crucial for the renormalization. In
section~5, we introduce the renormalization and proof estimates on the
N-block process. In section~6, we finally give the proof of
theorem~\procref{\prin}.
\subhead{\headnum Preliminaries}\endsubhead
In this section we introduce the notation used and the basic definitions.\headlabel{\pre}
\subsubhead{Norm and the lattice}\endsubsubhead We will use the $\Cal
L^1-$norm on $\Z^2$, that is, $|x-y|=\sum_{i=1,2}|x_i-y_i|$ for any
$x,y$ in $\Z^2$. For every subset $A$ of $\Z^2$ and $i=1,2$ we define
$\diam_i(A)=\sup\{|x_i-y_i|:x,\,y\in A\}$ and the diameter of $A$ is
$\diam(A)=\max(\diam_1(A), \diam_2(A))$. We turn $\Z^2$ into a graph
$(\Z^2,\E^2)$ with vertex set $\Z^2$ and edge set $\E^2=\{\{x, y\};
|x-y|=1\}$. If $x$ and $y$ are nearest neighbors, we denote this relation by $x\sim y$.
\subsubhead{Geometric objects}\endsubsubhead
A {\it box} $\Lambda$ is a finite subset of $\Z^2$ of the form
 $\Z^2\cap[a,b]\times[c,d]$. For $\und{r}\in (0,\infty)^2$, we define
 a box centered at the origin by
 $B(\und{r})=\Z^2\cap\Pi_{i=1,2}(-r_i/2, r_i/2]$. We say that the box
 is symmetric, if $r_1=r_2=r,$ and we denote it by $B(r)$. For
 $t\in\R^+$, we note the set ${\Cal H}_2(t)=\{\und{r}\in\R^2:
 r_i\in[t, 2t], i=1,2\}$. The set of all boxes in $\Z^2$, which are congruent to a box $B(\und{r})$
 with $\und{r}\in {\Cal H}_2(t)$, will be denoted by $\Cal B_2(t)$.
\subsubhead{Discrete topology}\endsubsubhead
Let $A$ be a subset of $\Z^2$. We define two different boundaries:
{\parindent 1cm
\item{$\bullet$}the inner vertex boundary: $\partial A=\{x\in
A|\;\exists y\in A^c \hbox{ such that } y\sim x\}$;
\item{$\bullet$}the edge boundary: $\partial^{\hbox{\sevenrm
edge}}A=\{\{x, y\}\in\E^2|\; x\in A,y\in A^c\}$. \par}
\noindent For a box $\Lambda$ and for each $i=\pm 1,\pm 2$, we define
the $i$th face $\partial_i\Lambda$ of $\Lambda$ by
$\partial_i\Lambda=\{x\in\Lambda|\;x_i\hbox{ is maximal}\}$ for $i$
positive and $\partial_i\Lambda=\{x\in\Lambda|\; x_{|i|}\hbox{ is
minimal}\}$ for $i$ negative. A {\it path} $\gamma$ is a finite or
infinite sequence $x_1, x_2,...$ of distinct nearest neighbors.
{\obeylines\subsubhead{FK percolation}\endsubsubhead
{\it Edge configurations.}}
The basic probability space for the edge processes is given by
$\Omega=\{0,1\}^{\E^2}$; its elements are called {\it edge
configurations in} $\Z^2$. The natural projections are given by
$\pr_e:\omega\in\Omega\mapsto\omega(e)\in\{0,1\}$, where
$e\in\E^2$. An edge $e$ is called open in the configuration $\omega$
if $\pr_e(\omega)=1$, and closed otherwise.\par
For $E\subseteq\E^2$ with $E\not=\emptyset$, we write $\Omega(E)$ for
the set $\{0,1\}^{E}$; its elements are called {\it configurations in}
$E$. Note that there is a one-to-one correspondence between cylinder
sets and configurations on finite sets $E\subset\E^2$, which is given
by $\eta\in\Omega(E)\mapsto\{\eta\}:=\{\omega\in\Omega\;|\;\omega(e)=\eta(e)\hbox{
for every } e\in E\}$.
We will use the following convention: the set $\Omega$ is regarded as
a cylinder (set) corresponding to the ``empty configuration'' (with
the choice $E=\emptyset.$) We will sometimes identify cylinders with
the corresponding configuration. For $A\subset\Z^2$, we set $\E(A)=\{(x,y):\
x,y\in A,\ x\sim y\}$. And let $\Omega_A$
stand for the set of the configurations in $A : \{0, 1\}^{\E(A)}$
and $\Omega^A$ for the set of the configurations {\it outside} $A :
\{0, 1\}^{\E^2\setminus\E(A)}$. In general, for $A\subseteq
B\subseteq\Z^2$, we set $\Omega_B^A=\{0,
1\}^{\E(B)\setminus\E(A)}$. Given $\omega\in\Omega$ and
$E\in\E^2$, we denote by $\omega(E)$ the restriction of $\omega$ to
$\Omega(E)$. Analogously, $\omega_B^A$ stands for the restriction of
$\omega$ to the set $\E(B)\setminus\E(A)$. \par Given
$\eta\in\Omega$, we denote by $\Cal O(\eta)$ the set of the edges of
$\E^2$ which are open in the configuration $\eta$. The connected
components of the graph $(\Z^2,\Cal O(\eta))$ are called
$\eta${\it-clusters}. The path $\gamma=(x_1,x_2,...)$ is said to be
$\eta$-open if all the edges $\{x_i, x_{i+1}\}$ belong to $\Cal
O(\eta)$. We write $\{A\leftrightarrow B\}$ for the event that there
exists an open path joining some site in $A$ with some site in
$B$.
\par If $V\subseteq\Z^2$ and $E$ consists of all the edges between
vertices in $V$, the graph $G=(V,E)\subseteq (\Z^2, \E^2)$ is called
the maximal subgraph of $(\Z^2, \E^2)$ on the vertices $V$.
Let $\omega$ be an edge configuration in $\Z^2$ (or in a subgraph of
$(\Z^2, \E^2)$). We can look at the {\it open clusters in} $V$ or
alternatively the open $V${\it-clusters}. These clusters are simply
the connected components of the random graph $(V, \Cal O
(\omega(E)))$, where $\omega(E)$ is th restriction of $\omega$ to
$E$.
\par For $A\subseteq B\subseteq\Z^2$, we use the notation $\Cal F_B^A$ for the
$\sigma$-field generated by the finite-dimensional cylinders
associated with configurations in $\Omega_B^A$. If $A=\emptyset$ or $B=\Z^2$, then we omit them from the notation.
\subsubhead{Stochastic domination}\endsubsubhead
There is a partial order $\preceq$ in $\Omega$ given by
 $\omega\preceq\omega^\prime$ iff $\omega(e)\leq\omega^{\prime}(e)$
 for every $e\in\E^2$. A function $f:\Omega\rightarrow\R$ is called
 {\it increasing} if $f(\omega)\leq f(\omega^{\prime})$ whenever
 $\omega\preceq\omega^{'}$. An event is called increasing if its
 characteristic function is increasing. Let $\Cal F$ be a
 $\sigma$-field of subsets of $\Omega$. For a pair of probability
 measures $\mu$ and $\nu$ on $(\Omega, \Cal F)$, we say that $\mu$
 {\it(stochastically) dominates} $\nu$ if for any $\Cal F$-measurable
 increasing function $f$ the expectations satisfy
 $\mu(f)\geq\nu(f)$.
\subsubhead{FK measures}\endsubsubhead
Let $V\subseteq\Z^2$ be finite and $E=\E(V).$ We first introduce {\it (partially
wired) boundary conditions} as follows. Consider a partition $\pi$ of the set $\partial
V$, say $\{B_1,...,B_n\}$. (The sets $B_i$) are disjoint nonempty subsets of $\partial V$
with $\bigcup_{i=1,...,n}B_i=\partial V$.) We say that $x, y\in\partial V$ are $\pi$-{\it
wired}, if $x,y\in B_i$ for an $i\in\{1,...,n\}$. Fix a configuration $\eta\in\Omega_V$.
We want to count the $\eta$-clusters in $V$ in such a way that $\pi$-wired sites are
considered to be connected. This can be done in the following formal way. We introduce an
equivalence relation on $V$: $x$ and $y$ are said to be $\pi\cdot\eta${\it-wired} if they
are both joined by $\eta$-open paths to (or identical with) sites $x^{\prime},
y^{\prime}\in\partial V$ which are themselves $\pi$-wired. The new equivalence classes
are called $\pi\cdot\eta$-clusters, or $\eta$-clusters in $V$ with respect to the
boundary condition $\pi$. The number of $\eta$-clusters in $V$ with respect to the
boundary condition $\pi$ (i.e., the number of $\pi\cdot\eta$-clusters) is denoted by
$\text{cl}^\pi(\eta)$. (Note that $\text{cl}^\pi$ is simply a random variable). For fixed
$p\in[0,1]$ and $q\geq 1$, the {\it FK measure on the finite set $V\subset\Z^2$ with parameters $(p,q)$ and boundary
conditions $\pi$} is a probability measure on the $\sigma$-field $\Cal F_V$, defined by
the formula
 $$
\forall \eta\in\Omega_V\qquad
\Phi_V^{\pi, p, q}[\{\eta\}]=\frac{1}{Z_V^{\pi, p,
q}}\left(\prod_{e\in
E}p^{\eta(e)}(1-p)^{1-\eta(e)}\right)q^{\text{cl}^{\pi}(\eta)},
\eqnum
 $$
 \eqlabel{\FKmeasure}\noindent
where $Z^{\pi, p, q}_V$ is the appropriate normalization factor.
Since $\Cal F_V$ is an atomic $\sigma$-field with atoms
 $\{\eta\},\eta\in\Omega_V$, \eqref{\FKmeasure} determines a unique
measure on $\Cal F_V$. Note that every cylinder has nonzero
probability.
There are two extremal b.c.s: the free boundary condition corresponds
to the partition $f$ defined to have exactly $|\partial V|$ classes,
and the wired b.c corresponds to the partition $w$ with only one
class. The set of all such measures called FK (or random cluster)
measures corresponding to different b.c.s will be denoted by
$\Cal R(p,q, V)$.
The stochastic process $(\pr_e)_{e\in\E(V)}:\Omega\rightarrow\Omega_V$ given on the
probability space $(\Omega, \Cal F, \Phi_V^{\pi, p, q})$ is called FK
{\it percolation with boundary conditions} $\pi$.
We will list some useful properties of FK measures with different
b.c.s. There is a partial order on the set of partitions of $\partial
V$. We say that $\pi$ {\it dominates} $\pi^{\prime},
\pi\geq\pi^{\prime},$ if $x,y$ $\pi^{\prime}$-wired implies that they
are $\pi$-wired. We then have $\Phi_V^{\pi^{\prime}, p, q}\preceq
\Phi_V^{\pi, p, q}$.
This implies immediately that for each $\Phi\in\Cal R(p, q, V),\, \Phi_V^{f, p,
q}\preceq\Phi \preceq\Phi_V^{w, p, q}.$
Next we discuss properties of conditional FK measures.
For given $U\subseteq V$ and $\omega\in\Omega$, we define a partition
$W_V^U(\omega)$ of $\partial U$ by declaring $x, y\in\partial U$ to be
$W_V^U(\omega)$-wired if they are joined by an $\omega_V^U$-open path.
Fix a partition $\pi$ of $\partial V$. We define a new partition of $\partial U$ to
be $\pi\cdot W_V^U(\omega)$-wired if they are both joined by $\omega_V^U$-open paths to
(or identical with) sites $x^{\prime}, y^{\prime}$, which are themselves $\pi$-wired.
Then, for every $\Cal F_U$-measurable function $f$,
$$
\Phi_V^{\pi, p,q}[f | \Cal F_V^U](\omega)=\Phi_V^{\pi\cdot
W_V^U(\omega), p,q}[f],\qquad
\Phi_V^{\pi, p,q} \hbox{ a.s.}
\eqnum $$ \eqlabel{\FKmeasuretwo}
Note that \eqref{\FKmeasuretwo} can be interpreted as a kind of Markov
property.
A direct consequence is the {\it finite-energy property}. Fix an edge $e$ of
$\E(V)$ and denote by $\Cal F_V^e$ the $\sigma$-algebra generated by the
random variables $\{pr_b;b\in\E(V)\setminus\{e\}\}$. Then
$$
\eqalign{
\Phi&_V^{\pi,p,q}[e\hbox{ is open }|\Cal F_V^e](\omega)
=\left\{\eqalign{&p\hbox{ if the endpoints of }e\hbox{ are }\pi\cdot W^e_V\hbox{-wired},\cr
&p/[p+q(1-p)]\hbox{ otherwise}.}
\right.}\eqnum$$\eqlabel{\NRJ}\noindent
The equality \eqref{\FKmeasuretwo} leads to volume monotony for
FK-measures. Let $U\subset V$, for every increasing function $g\in\Cal F_U$ and
$\Phi_V\in\Cal R(p,q,V),$ we have that
$$\eqalign{
&\Phi_U^{f,p,q}[g]\leq\Phi_V[g\;|\;\Cal F_V^U]\leq\Phi_U^{w,p,q}[g]\quad\Phi_V\hbox{ a.s. ,}\cr
&\Phi_U^{f,p,q}[g]\leq\Phi_V^{f,p,q}[g]\leq\Phi_V^{w,p,q}[g]\leq\Phi_U^{w,p,q}[g].
}$$
\subsubhead{Planar duality for FK-measures}\endsubsubhead
Let $n>0$ and $\und{n}\in\Cal H_2(n)$. To the set
 $B(\und{n})\subset\Z^2$ we associate the set
 $\widehat{B}(\und{n})\subset\Z^2+(1/2,1/2)$, which is defined as the
 smallest box of $\Z^2+(1/2,1/2)$ containing $B(\und{n})$, see the
 figure 1 below.
Notice that if $\und{n}\in\Cal H_2(n)$ then
 $\widehat{B}(\und{n})\in[\Cal B_2(n+1)\cup\Cal B_2(n)]+(1/2,1/2)$. To each edge $e\in\E(B(\und{n}))$ we associate the edge
 $\widehat{e}\in\E(\widehat{B}(\und{n}))$ that crosses the edge
 $e$. According to \cite{\GrimoII}, if we associate to each
 configuration $\omega\in\Omega_{B(\und{n})}$ the dual configuration $\widehat{\omega}_\omega$:
 $$
\widehat{\omega}_\omega\in\Omega_{\widehat{B}(\und{n})}^{\partial\widehat{B}(\und{n})}\,\hbox{
  such that }
\forall e\in\E(B(\und{n})),\ \widehat{\omega}_\omega(\widehat{e})=1-\omega(e),
$$
then we have that
$$
\Phi^{f,p,q}_{B(\und{n})}[\omega]=\Phi^{w,\widehat{p},q}_{\widehat{B}(\und{n})}[\{\omega_d\in\Omega_{\widehat{B}(\und{n})},\
\forall\widehat{e}\in\E(\widehat{B}(\und{n}))\setminus\E(\partial\widehat{B}(\und{n})):\
\omega_d(\widehat{e})=\widehat{\omega}_\omega(\widehat{e})\}],
$$
where $\widehat{p}$ is the dual point of $p$ : $\widehat{p}=q(1-p)/[p+q(1-p)]$.\par
\vbox{
\centerline{
\hbox{
\psset{unit=0.5cm}
\pspicture(-8,-8)(8,8)
\psline[linestyle=dashed]{-}(-6,-5)(6,-5)
\psline[linestyle=dashed]{-}(-6,-3)(6,-3)
\psline[linestyle=dashed]{-}(-6,-1)(6,-1)
\psline[linestyle=dashed]{-}(-6,1)(6,1)
\psline[linestyle=dashed]{-}(-6,3)(6,3)
\psline[linestyle=dashed]{-}(-6,5)(6,5)
\psline[linestyle=dashed]{-}(-6,-5)(-6,5)
\psline[linestyle=dashed]{-}(-4,-5)(-4,5)
\psline[linestyle=dashed]{-}(-2,-5)(-2,5)
\psline[linestyle=dashed]{-}(0,-5)(0,5)
\psline[linestyle=dashed]{-}(2,-5)(2,5)
\psline[linestyle=dashed]{-}(4,-5)(4,5)
\psline[linestyle=dashed]{-}(6,-5)(6,5)
\psdots[dotstyle=o](-6,-5)(-4,-5)(-2,-5)(0,-5)(2,-5)(4,-5)(6,-5)
\psdots[dotstyle=o](-6,-3)(-4,-3)(-2,-3)(0,-3)(2,-3)(4,-3)(6,-3)
\psdots[dotstyle=o](-6,-1)(-4,-1)(-2,-1)(0,-1)(2,-1)(4,-1)(6,-1)
\psdots[dotstyle=o](-6,1)(-4,1)(-2,1)(0,1)(2,1)(4,1)(6,1)
\psdots[dotstyle=o](-6,3)(-4,3)(-2,3)(0,3)(2,3)(4,3)(6,3)
\psdots[dotstyle=o](-6,5)(-4,5)(-2,5)(0,5)(2,5)(4,5)(6,5)
\psline{-}(-5,-4)(5,-4)
\psline{-}(-5,-2)(5,-2)
\psline{-}(-5,0)(5,0)
\psline{-}(-5,2)(5,2)
\psline{-}(-5,4)(5,4)
\psline{-}(-5,-4)(-5,4)
\psline{-}(-3,-4)(-3,4)
\psline{-}(-1,-4)(-1,4)
\psline{-}(1,-4)(1,4)
\psline{-}(3,-4)(3,4)
\psline{-}(5,-4)(5,4)
\psdots(-5,-4)(-3,-4)(-1,-4)(1,-4)(3,-4)(5,-4)
\psdots(-5,-2)(-3,-2)(-1,-2)(1,-2)(3,-2)(5,-2)
\psdots(-5,0)(-3,0)(-1,0)(1,0)(3,0)(5,0)
\psdots(-5,2)(-3,2)(-1,2)(1,2)(3,2)(5,2)
\psdots(-5,4)(-3,4)(-1,4)(1,4)(3,4)(5,4)
\uput[-180](-3.5,-6.1){$e\in\E(B(\und{n}))$}
\psline[linewidth=0.5pt]{->}(-3.75,-5.9)(-2.5,-4)
\uput[0](0.5,-6){$\widehat{e}\in \E(\widehat{B}(\und{n}))\setminus\E(\partial\widehat{B}(\und{n}))$}
\psline[linewidth=0.5pt]{->}(0.5,-5.9)(-2,-4.5)
\endpspicture
}
}
\centerline{figure 1: A box and its dual}
}\medskip\noindent

\par\noindent Thus, for each $\Cal F_{B(\und{n})}$ measurable event
$A$ we can associate a $\Cal
 F_{\widehat{B}(\und{n})}^{\partial\widehat{B}(\und{n})}$
 measurable event
 $$
\widehat{A}=\{\omega_d\in\Omega_{\widehat{B}(\und{n})}:\
 \exists\ \omega\in A,
\forall\widehat{e}\in\E(\widehat{B}(\und{n}))\setminus\E(\partial\widehat{B}(\und{n})):\
\omega_d(\widehat{e})=\widehat{\omega}_\omega(\widehat{e})\},
$$
which satisfies
$$
\Phi^{f,p,q}_{B(\und{n})}[A]=\Phi^{w,\widehat{p},q}_{\widehat{B}(\und{n})}[\widehat{A}].
$$
 \subhead{\headnum  Connectivity in boxes}\endsubhead
In this section we establish preliminary estimates on crossing events in boxes. We rely on the exponential decay of the connectivities in the dual subcritical model.\headlabel{\connect}
The usual definition of the exponential decay is based on the infinite
volume FK-measure $\Phii^{p,q}$. But we are concerned by asymptotics of finite
volume measures and we would like to use the exponential decay in
finite boxes.
In order to translate the exponential decay to the finite volume
measures we need a control on
the effects of boundary conditions. As shown in \cite{\AlexI}, the infinite FK-measure on $\Z^2$ satisfies the weak mixing
property as soon as the connectivities decay exponentially. That is to say for
all
events $A, B$ which are respectively $\Cal F_\Lambda$ measurable and $\Cal F_\Gamma$
measurable with $\Lambda, \Gamma\subseteq \Z^2$ then $|\Phii^{p,q}[A\cap
B]-\Phii^{p,q}[A]\Phii^{p,q}[B]|$ decreases exponentially in the distance between
$\Lambda$ and $\Gamma$. This weak mixing property implies, as proved in
\cite{\AlexII}, that we have exponential decay in finite boxes as soon as the exponential
decay for the infinite volume measure holds ($p<\pg$):
\proclaim{Proposition \procnum\ (Theorem 1.2 of \cite{\AlexII})}
 Let $q\geq 1$ and $p<\pg$. There exists two positive constants $c$
 and $\lambda$ such that for all boxes $\Lambda\subset\Z^2$ and for all $x,
 y$ in $\Lambda$, we have that
 $$
 \Phi^{w,p,q}_\Lambda[x\leftrightarrow y\hbox{ in }\Lambda]\leq\lambda\exp(-c|x-y|).
 $$\endproclaim\noindent
In fact, theorem 1.2 of \cite{\AlexII} is more general and applies to sets $\Lambda$
which are not boxes and to general boundary conditions. From this result, we
get that
\proclaim{Lemma \procnum}
Let $q\geq 1$ and $p<\pg$. There  exists a positive constant $c$ such
that for all positive integers $n$ and for $l$ large enough, we have that
 $$
 \sup_{\und{n}\in \Cal H_2(n)}
\Phi^{w,p,q}_{B(\und{n})}[\exists\hbox{ an open path in }B(\und{n})\hbox{ of
  diameter }\geq l]
\leq n^2\exp(-c l).
 $$
\endproclaim
\proclabel{\Probdualpath}\noindent
\demo{Proof}
Let us fix $n$ and $l$, then we have
$$\eqalign{
&\sup_{\und{n}\in \Cal H_2(n)}\Phi^{w,p,q}_{B(\und{n})}[\exists\hbox{ an open path in }B(\und{n})\hbox{ of diameter }\geq l]\cr
\leq & 4n^2 \sup_{\und{n}\in \Cal H_2(n)}\sup_{x\in B(\und{n})}\Phi^{w,p,q}_{B(\und{n})}[x\leftrightarrow \partial B(x,2l)\hbox{ in }B(\und{n})]\cr
\leq & 32n^2l \sup_{\und{n}\in \Cal H_2(n)}\sup_{x\in
B(\und{n})}\sup_{y\in\partial B(x,2l)}\Phi^{w,p,q}_{B(\und{n})}[x\leftrightarrow y\hbox{ in
}B(\und{n})]\cr \leq & 32\lambda n^2l \exp(-c l),
}
$$
where we used proposition 2 in the last line. The result follows
by taking $l$ large enough.
\qed
\enddemo
\noindent As a first consequence of the exponential decay in finite boxes, we obtain:
\proclaim{Lemma \procnum}
For $p>\pgd$ we have,
$$
\lim_{n\rightarrow\infty}\Phi^{f, p,q}_{B(n)}[0\leftrightarrow\partial B(n)]=\theta(p,q).
$$\proclabel{\thdeux}\endproclaim
\demo{Proof}
Let $N<n$, then
 $$
\eqalign{
\Phi^{f, p, q}_{B(n)}[0\leftrightarrow
\partial B(N)]-\Phi^{f,p,q}_{B(n)}[0\leftrightarrow \partial
B(N)\,,\;0\nleftrightarrow\partial B(n)]
 =&\Phi^{f, p, q}_{B(n)}[0\leftrightarrow\partial B(n)]\cr
 \leq&\Phi^{f, p,q}_{B(n)}[0\leftrightarrow \partial B(N)]. }
\eqnum
$$
\eqlabel{\sarko}\noindent
 Now we will estimate $\Phi^{f, p,q}_{B(n)}[0\leftrightarrow \partial B(N)\,,\;
0\nleftrightarrow \partial B(n)]$: by symmetry, $$
\Phi^{f, p,q}_{B(n)}[0\leftrightarrow \partial B(N)\,,\;
0\nleftrightarrow\partial B(n)]\leq 4 \Phi^{f, p,q}_{B(n)}[0\leftrightarrow \partial_1 B(N)\,,\; 0\nleftrightarrow
\partial B(n)].
$$
\noindent Then for $N$ large enough we have that
$$\eqalign{
\Phi^{f, p, q}_{B(n)}[0\leftrightarrow \partial_1
B(N),0\nleftrightarrow\partial B(n)]\leq&\Phi^{f, p, q}_{B(n)}
\left[
\eqalign{
&\exists k>0\ \exists j\geq N : \exists\hbox{ an open dual path}\cr
&\hbox{joining
}(-k+\frac{1}{2},\frac{1}{2})\hbox{ to }(j+\frac{1}{2},\frac{1}{2}) }\right]\cr
\leq&\sum_{k>0,\, j\geq N}\exp(-c(k+j))\cr
\leq&\exp(-cN),
 }$$
for a certain positive constant $c$.
\par\noindent By taking the limit $n\rightarrow\infty$ in \eqref{\sarko} we get
 $$
 \eqalign{
 \Phi^{p, q}_{\infty}[0\leftrightarrow \partial B(N)]-4e^{-dN}\leq &
\liminf_{n\rightarrow\infty}\Phi^{f, p, q}_{B(n)}[0\leftrightarrow\partial B(n)]\cr \leq
& \limsup_{n\rightarrow\infty}\Phi^{f, p,q}_{B(n)}[0\leftrightarrow\partial B(n)]
\leq  \Phi^{p, q}_{\infty}[0\leftrightarrow\partial B(N)], }
 $$
finally by taking the limit $N\rightarrow\infty$, we get the desired result.
\qed
\enddemo
Next, we define events that will be crucial in the renormalization procedure. For this,
we introduce the notion of {\it crossing}. Let $B\subset\Z^2$ be a finite box. For
$i=1,2$ we say that a $i$-{\it crossing} occurs in $B$, if $\partial_{-i}B$ and
$\partial_{i}B$ are joined by an open path in $B$. In addition to that, we say that a
cluster $C$ of $B$ is crossing in $B$, if $C$ contains a $1$-crossing path and a
$2$-crossing path.
\par\noindent For $\und{n}\in{\Cal H}_2(n)$, we set
 $$
U(\und{n})=\{\exists! \hbox{ open cluster  }C^* \hbox{ crossing } B(\und{n})\}.
 $$
For a monotone, increasing function $g: \N\tend[0,\infty)$ with $g(n)\leq n$, let us
define
 $$
R^g(\und{n})=U(\und{n})\cap
\left\{\eqalign{
&\hbox{ every open path } \gamma\subset B(\und{n})\hbox{ with }\cr
&\diam(\gamma)\geq g(n)\hbox{ is contained in } C^* } \right\}.
 $$
And finally we set
 $$
O^g(\und{n})=R^g(\und{n})\cap
\left\{\eqalign{ &C^*\hbox{ crosses every sub-box }\cr
& \Cal Q \in \Cal B_2(g(n))\hbox{ contained in }B(\und{n})}\right\}.
 $$
The next theorem gives the desired estimates on the above mentioned events.
\proclaim{Theorem \procnum} Assume $p>\pgd$. We have
 $$
\limsup_{n\rightarrow \infty} \frac{1}{n}\log\sup_{\und{n}\in \Cal H_2(n)}
\sup_{\Phi\in \Cal R(p,q, B(\und{n}))}\Phi[U(\und{n})^c]<0.
\eqnum$$
\eqlabel{\devU}
Also, there exists a constant $\kappa=\kappa(p,q)>0$ such that $\liminf_{n\tend\infty} g(n)/\log n
>\kappa$ implies $$
\limsup_{n\rightarrow \infty}\frac{1}{g(n)}\log\sup_{n\in\Cal H_2(n)} \sup_{\Phi\in\Cal R(p,q,B(\und{n}))}\Phi[R^g(\und{n})^c]<0.\eqnum
$$
\eqlabel{\devR}
There exists a constant $\kappa'=\kappa'(p,q)>0$ such that $\liminf_{n\tend\infty} g(n)/\log n>\kappa'$
implies
 $$
 \limsup_{n\rightarrow \infty} \frac{1}{g(n)}\log\sup_{n\in\Cal H_2(n)}
\sup_{\Phi\in\Cal R(p,q,B(\und{n}))}\Phi[O^g(\und{n})^c]<0.
\eqnum
$$
\eqlabel{\devO}
\endproclaim \proclabel{\thun}
\proclaim{Remark}
Note that in dimension two, if there is a crossing cluster then it is unique.
\endproclaim
\demo{Proof}
As $U(\und{n})^c$ is decreasing we have for every $\Phi\in\Cal R(p,q, B(\und{n}))$ that
 $$\eqalign{
\Phi[U(\und{n})^c]\leq&\Phi^{f,p,q}_{B(\und{n})}[U(\und{n})^c]\cr
\leq&\Phi^{f,p,q}_{B(\und{n})}[\nexists\;1\hbox{-crossing for
}B(\und{n})]+\Phi^{f,p,q}_{B(\und{n})}[\nexists\; 2\hbox{-crossing for
}B(\und{n})]\cr
\leq&\sum_{i=1,2}\Phi^{f,p,q}_{B(\und{n})}[\partial_{-i}\widehat{B}(\und{n})\leftrightarrow\partial_{i}\widehat{B}(\und{n})\hbox{
in
}\widehat{B}(\und{n})\setminus\partial\widehat{B}(\und{n})],
 }$$
the last inequality follows from planar duality: if there is no
$1$-crossing in the original lattice then  $\partial_{-2}\widehat{B}(\und{n})\leftrightarrow\partial_{2}\widehat{B}(\und{n})\hbox{
in }\widehat{B}(\und{n})\setminus\partial\widehat{B}(\und{n})$ for the
corresponding dual configuration. The same argument works for the
$2$-crossing. Thus, we have that
$$
\Phi[U(\und{n})^c]\leq 2\Phi^{w,\widehat{p},q}_{\widehat{B}(\und{n})}[\exists\hbox{ an open path in }\widehat{B}(\und{n})\hbox{ of diameter }\geq n],
$$
and \eqref{\devU} follows from lemma \procref{\Probdualpath}.
\par\noindent For the second inequality, let us note that $$
 R^g(\und{n})^c\subset U(\und{n})^c\bigcup\left(U(\und{n})\cap\left\{\eqalign{
 &\exists\hbox{ an open path }\gamma\hbox{ of }B(\und{n})\hbox{ with }\cr
 &\diam(\gamma)\geq g(n)\hbox{ not contained in } C^*}\right\}\right).
$$
By \eqref{\devU}, we have only to consider the second term.
\par\noindent
In
 $$U(\und{n})\cap\left\{\eqalign{
 &\exists\hbox{ an open path }\gamma\hbox{ of }B(\und{n})\hbox{ with }\cr
 &\diam(\gamma)\geq g(n)\hbox{ not contained in } C^*}\right\},
 $$
by proposition 11.2 of \cite{\Grimo} and by considering all the edges of
$\E(\partial \widehat{B}(\und{n}))$ open, there is a unique innermost open dual
circuit containing $\gamma$ in its
interior. From this dual circuit, we extract an open dual path living
in the graph
$(\widehat{B}(\und{n}),\E(\widehat{B}(\und{n}))\setminus\E(\partial\widehat{B}(\und{n})))$
of diameter greater than $g(n)$: Without lost of generality, we can
suppose that $\diam(\gamma)=\diam_1(\gamma)$ and that
$\gamma\nleftrightarrow\partial_2 B(\und{n})$. Among the vertices of the dual
circuit surrounding $\gamma$, let  $\widehat{x}$ be the highest vertex among the most
on the left, and let $\widehat{y}$ be the highest vertex among the most on the
right. Then there is an arc joining $\widehat{x}$ and $\widehat{y}$ in
$(\widehat{B}(\und{n}),\E(\widehat{B}(\und{n}))\setminus\E(\partial\widehat{B}(\und{n})))$.
This arc is of diameter larger than $g(n)$.
 Thus by lemma \procref{\Probdualpath} there is a positive constant $c$ such that for
$n$ large enough we have that
$$
\Phi\left[U(\und{n})\cap
\left\{\eqalign{
 &\exists\hbox{ an open path }\gamma\hbox{ of }B(\und{n})\hbox{ with }\cr
 &\diam(\gamma)\geq g(n)\hbox{ not contained in } C^*}\right\}\right]
\leq n^2\exp[-cg(n)].
$$ Now, take an $\alpha>0$ such that $\alpha c>1$, then for $g$ such that
$g(n)>2\alpha\log n /(\alpha c-1)$ we have that $$
\limsup_{n\tend\infty}\frac{1}{g(n)}\log(n^2\exp[-cg(n)])<-\frac{1}{\alpha},
$$
which concludes the proof of \eqref{\devR}.
\par\noindent To study $O^g(\und{n})$, we remark that the number of boxes
$\Cal Q$ of $\Cal B_2(g(n))$ contained
 in $B(\und{n})$ is bounded by $16n^4$. This implies that for every $\Phi\in \Cal R(p,q,B(\und{n}))$ one gets
 $$\eqalign{
\Phi[O^g(\und{n})^c]\leq& \Phi[R^g(\und{n})^c] + 16n^4\sup_{\Cal{Q}\in
\Cal{B}_2(g(n))} \Phi[\nexists \hbox{ crossing in }\Cal{Q}]\cr
\leq&\Phi[R^g(\und{n})^c] + 16n^4\sup_{\Cal{Q}\in
\Cal{B}_2(g(n))} \Phi^{f,p,q}_{B(\und{n})}[\nexists \hbox{ crossing in }\Cal{Q}]\cr
\leq&\Phi[R^g(\und{n})^c] + 16n^4\; \sup_{\Cal{Q}\in
\Cal{B}_2(g(n))}\Phi^{f,p,q}_{\Cal{Q}}[\nexists \hbox{ crossing in }\Cal{Q}].
 }$$
To deduce the last inequality, we notice that $\{\nexists \hbox{ crossing in }\Cal{Q}\}$
is a decreasing event and that all the $\Cal{Q}\in\Cal{B}_2(g(n))$ are smaller than
$B(n)$, thus for all $\Cal Q\in\Cal B_2(g(n))$ that are included in $B(n)$ we
have that
 $$
\Phi^{f,p,q}_{B(\und{n})}[\nexists \hbox{ crossing in }\Cal{Q}]\leq
\Phi^{f,p,q}_{\Cal{Q}}[\nexists \hbox{ crossing in }\Cal{Q}].
$$
\par\noindent The first term in the r.h.s. has been treated previously. By
\eqref{\devU} the second term is bounded by $n^4\exp[-c g(n)]$ for a certain positive
constant $c$ and we conclude the proof as before.\qed
\enddemo
\subhead{\headnum  Renormalization}\endsubhead
In this section we adapt the renormalization procedure introduced in \cite{\Pisz} to the
two dimensional case. For this, let $N\geq 24$ be an integer.\headlabel{\renorma} We say that a subset
$\Lambda$ of $\Z^2$ is a $N$-{\it large box} if $\Lambda$ is a
finite box containing a symmetric box of scale-length $3N$, i.e., if
$\Lambda=\Z^2\cap\prod_{i=1,2}(a_i,b_i]$ where $b_i-a_i\geq 3N$ for $i=1,2.$ When
$\Lambda$ is a $N$-large box, one can partition it with blocks of $\Cal B(N)$. We
first define the $N$-{\it rescaled box} of $\Lambda$: $\Lambda^{(N)}=\{\bl{k}
\in\Z^2\;|\;T_{N\bl{k}}(-N/2,N/2]^2\subseteq\Lambda\};$ where $T_a$ is the translation in
$\Z^2$ by a vector $a\in\Z^2$. We turn $\Lambda^{(N)}$ into a graph by endowing it with
the set of edges $\E(\Lambda^{(N)})$. Then we define the partitioning  blocks:
\item{$\bullet$} If $\bl{k}\in\Lambda^{(N)}\setminus\partial\Lambda^{(N)}$ then
$B_\bl{k}=T_{N\bl{k}}(-N/2,N/2]^2$.\par
\item{$\bullet$} If $\bl{k}\in\partial\Lambda^{(N)}$ then some care is needed in order to get
a partition. In this case we define the set
 $$
\Cal
M(\bl{k})=\{\bl{l}\in\Z^2\;|\;\bl{l}\sim\bl{k},
T_{N\bl{l}}(-N/2,N/2]^2\cap\Lambda\not=\emptyset,\;T_{N\bl{l}}(-N/2,N/2]^2\cap\Lambda^{c}\not=\emptyset\},
 $$
and the corresponding blocks become
 $$
B_\bl{k}=T_{N\bl{k}}(-N/2,N/2]^2\cup\bigcup_{\bl{l}\in\Cal
M(\bl{k})}\left(T_{N\bl{l}}(-N/2,N/2]^2\cap\Lambda\right).
 $$
\par\noindent The collection of sets $\{B_\bl{k},\ \bl{k}\in\Lambda^{(N)}\}$ is a partition of $\Lambda$
into blocks included in $\Cal B (N)$, see figure 2.\par
\midinsert
\centerline{
\psset{unit=0.3cm}
\pspicture(-12,-10)(12,12)
\psline{-}(-9.6,-4)(9.6,-4)
\psline{-}(-9.6,-2)(9.6,-2)
\psline{-}(-9.6,+0)(9.6,+0)
\psline{-}(-9.6,+2)(9.6,+2)
\psline{-}(-9.6,+4)(9.6,+4)
\psline{-}(-7,-7.6)(-7,7.6)
\psline{-}(-5,-7.6)(-5,7.6)
\psline{-}(-3,-7.6)(-3,7.6)
\psline{-}(-1,-7.6)(-1,7.6)
\psline{-}(+1,-7.6)(+1,7.6)
\psline{-}(+3,-7.6)(+3,7.6)
\psline{-}(+5,-7.6)(+5,7.6)
\psline{-}(+7,-7.6)(+7,7.6)
\uput[-90](0,10){$\Lambda$}
\pspolygon[linewidth=2pt](-9.6,-7.6)(-9.6,7.6)(9.6,7.6)(9.6,-7.6)(-9.6,-7.6)
\psset{dotsize=0.4}
\psline[linewidth=0.5pt]{->}(-11,1.5)(-8,-1)
\uput[-180](-11,1.5){$\bold{k}\in\Lambda^{(N)}$}
\psline[linewidth=0.5pt]{->}(-11,-2)(-8.5,-1.5)
\uput[-180](-11,-2){$B_\bold{k}$}
\psdots(-8,-5)(-6,-5)(-4,-5)(-2,-5)(-0,-5)(2,-5)(4,-5)
(6,-5)(8,-5)
\psdots(-8,-3)(-6,-3)(-4,-3)(-2,-3)(-0,-3)(2,-3)(4,-3)
(6,-3)(8,-3)
\psdots(-8,-1)(-6,-1)(-4,-1)(-2,-1)(-0,-1)(2,-1)(4,-1)
(6,-1)(8,-1)
\psdots(-8,1)(-6,1)(-4,1)(-2,1)(-0,1)(2,1)(4,1)
(6,1)(8,1)
\psdots(-8,3)(-6,3)(-4,3)(-2,3)(-0,3)(2,3)(4,3)
(6,3)(8,3)
\psdots(-8,5)(-6,5)(-4,5)(-2,5)(-0,5)(2,5)(4,5)
(6,5)(8,5)
\endpspicture
}
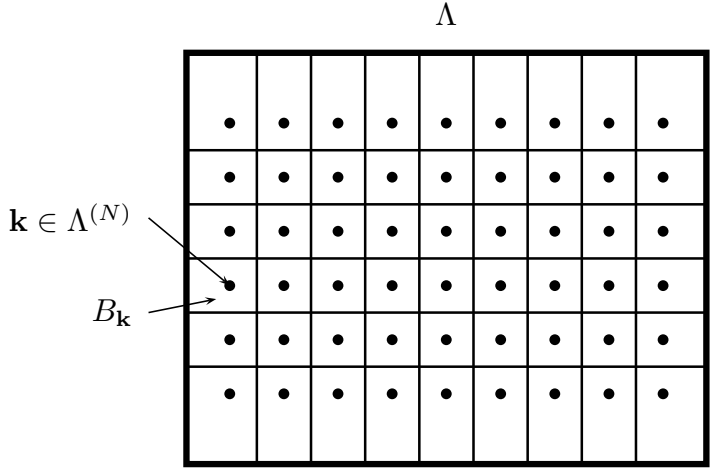
\centerline{figure 2: The partition of $\Lambda$}
\endinsert\medskip\noindent
In addition to the boxes $\{B_\bl{k},\ \bl{k}\in\Lambda^{(N)}\}$ we associate to each edge $(\bl{k},\bl{l})$
of $\E(\Lambda^{(N)})$ the box $D_{(\bl{k}, \bl{l})}$. More precisely, for $(\bl{k},\bl{l})\in\E(\Lambda^{(N)})$ such that $\sum_{j=1,2}|\bl{k}_j-\bl{l}_j|=\bl{k}_i-\bl{l}_i=1$, we
define $m(\bl{l},\bl{k})=T_{N\bl{l}}(\lfloor N/2\rfloor e^{(i)})$, where $(e^{(1)},
e^{(2)})$ is the canonical orthonormal base of $\Z^2$ and $\lfloor r
\rfloor$ denotes the integer part of $r$. The point $m(\bl{l},\bl{k})$ represents the
middle of the $i$-th face of $B_\bl{l}$. We then define the box
$D_{(\bl{l},\bl{k})}=D_{(\bl{k},\bl{l})}=T_{m(\bl{l},\bl{k})}(B(\lfloor N/4\rfloor))$.
\par Now we have all the needed geometric objects to construct our
renormalized (dependent) site percolation process on $(\Lambda^{(N)},\E(\Lambda^{(N)}))$.
This process will depend on the original FK-percolation process only through a number of
events defined in the boxes $(B_\bl{k})_{\bl{k}\in\Lambda^{(N)}}$ and
$(D_{e})_{e\in\E(\Lambda^{(N)})}$. These events are:
\item{$\bullet$} For all $(\bl{k}, \bl{l})\in \E(\Lambda^{(N)})$ such that
$\sum_{j=1,2}|\bl{k}_j-\bl{l}_j|=\bl{k}_i-\bl{l}_i=1$, we define
$$
 K_{\bl{k},\bl{l}}=\{\exists\ i\hbox{-crossing in }
D_{\bl{k}, \bl{l}}\},\qquad
K_\bl{k}=\bigcap_{\bl{j}\in\Lambda^{(N)}:\bl{j}\sim\bl{k}}K_{\bl{k}, \bl{j}}.
$$ \par
\item{$\bullet$} For all $\bl{i}\in\Lambda^{(N)}$, we define
$$
 R_\bl{i}=\{\exists! \hbox{ a
crossing cluster } C^*_\bl{i}\hbox{ in }B_\bl{i}\}\cap
 \left\{\eqalign{
 &\hbox{every open path } \gamma\subset B_\bl{i}\hbox{ with }\cr
 &\diam(\gamma)\geq\frac{\sqrt{N}}{10}\hbox{ is included in }C^*_{\bl{i}}
 }\right\}.
$$
\par Finally our renormalized process is  the indicator of the occurrence of the above mentioned
events:
$$
\forall\bl{k}\in\Lambda^{(N)}\quad X_\bl{k}=\left\{\matrix
 1 & \hbox{ on }\hfill & R_\bl{k}\cap K_\bl{k}\cr
0 & \hbox{otherwise}\hfill\cr \endmatrix\right.
 $$
We also call the process $\{X_\bl{k},\bl{k}\in\Lambda^{(N)}\}$ the $N$-block
process and whenever $X_\bl{k}=1$, we say that the block  $B_\bl{k}$ is occupied. As explained in \cite{\Pisz}, the  $N$-block process has the following important
geometrical property: if $C^{(N)}$ is a cluster of occupied blocks then there is a
unique cluster $C$ of the underlying microscopic FK-percolation process that crosses
all the blocks $\{B_\bl{k},\ \bl{k}\in C^{(N)}\}$. Moreover, the events involved in the definition of the $N$-block process become more probable as the size of the
blocks increases. This leads us to the following stochastic domination result:
\proclaim{Proposition \procnum}
Let $q\geq 1$ and $p>\pgd$. Then for $N$ large enough, every $N$-large  box $\Lambda$ and every measure $\Phi^\pi\in\Cal R(p, q, \Lambda)$, the law of the
$N$-block process $(X_\bl{i})_{\bl{i}\in\Lambda^{(N)}}$ under $\Phi^\pi$,
stochastically dominates independent site percolation on $\Lambda^{(N)}$ with parameter
$p(N)=1-\exp(-C \sqrt{N})$, where $C$ is a positive constant. \endproclaim
\proclabel{\SM}
\demo{Proof}
According to \cite{\LSS}, it is sufficient to establish that for $N$
large enough and for all $\bl{i}\in\Lambda^{(N)}$ the following
inequality holds:
 $$
\Phi^{\pi}[X_\bl{i}=0\;|\;\sigma(X_\bl{j}: |\bl{j}-\bl{i}|>1)]\leq \exp(-C \sqrt{N}). \eqnum
 $$ \eqlabel{\ineqcond}\par
\noindent In what follows, we use the same notation for positive constants that
may differ from one line to another.
In order to prove \eqref{\ineqcond}, we will consider the set $$\Dsty
E_\bl{i}=\bigcup_{|\bl{j}-\bl{i}|\leq 1}B_\bl{j}
 \setminus\bigcup_{\bl{j}\sim\bl{i}}\;\bigcup_{\bl{k}\sim\bl{j},
k\not=\bl{i}}D_{\bl{j},\bl{k}},$$
as drawn in figure 3.\par
\midinsert
\centerline{
\hbox{
\psset{unit=0.25cm}
\pspicture(-12,-13)(12,14)
\psline[linewidth=1.5pt]{-}(-4,-12)(-1,-12)(-1,-11)(1,-11)(1,-12)(4,-12)
\psline[linewidth=1.5pt]{-}(4,-12)(4,-9)(3,-9)(3,-7)(4,-7)(4,-4)
\psline[linewidth=1.5pt]{-}(4,-4)(7,-4)(7,-3)(9,-3)(9,-4)(12,-4)
\psline[linewidth=1.5pt]{-}(12,-4)(12,-1)(11,-1)(11,1)(12,1)(12,4)
\psline[linewidth=1.5pt]{-}(12,4)(9,4)(9,3)(7,3)(7,4)(4,4)
\psline[linewidth=1.5pt]{-}(4,4)(4,7)(3,7)(3,9)(4,9)(4,12)
\psline[linewidth=1.5pt]{-}(4,12)(1,12)(1,11)(-1,11)(-1,12)(-4,12)
\psline[linewidth=1.5pt]{-}(-4,12)(-4,9)(-3,9)(-3,7)(-4,7)(-4,4)
\psline[linewidth=1.5pt]{-}(-4,4)(-7,4)(-7,3)(-9,3)(-9,4)(-12,4)
\psline[linewidth=1.5pt]{-}(-12,4)(-12,1)(-11,1)(-11,-1)(-12,-1)(-12,-4)
\psline[linewidth=1.5pt]{-}(-12,-4)(-9,-4)(-9,-3)(-7,-3)(-7,-4)(-4,-4)
\psline[linewidth=1.5pt]{-}(-4,-4)(-4,-7)(-3,-7)(-3,-9)(-4,-9)(-4,-12)
\psline[linestyle=dashed,linewidth=1pt]{-}(-4,-4)(-4,4)(4,4)(4,-4)(-4,-4)
\psline[linewidth=1pt]{-}(-1,5)(1,5)(1,3)(-1,3)(-1,5)
\psline[linewidth=1pt]{-}(-1,-3)(1,-3)(1,-5)(-1,-5)(-1,-3)
\psline[linewidth=1pt]{-}(3,1)(5,1)(5,-1)(3,-1)(3,1)
\psline[linewidth=1pt]{-}(-5,1)(-3,1)(-3,-1)(-5,-1)(-5,1)
\psline[linewidth=0.5pt]{->}(-6.5,-8.5)(-2,-2)
\uput[-180](-6,-9){$B_{\bold{i}}$}
\psline[linewidth=0.5pt]{->}(8,-7)(4.5,0)
\uput[-90](8,-7){$D_{\bold{i},\bold{j}}$}
\psline[linewidth=0.5pt]{->}(-8,9)(-4,5.5)
\uput[-180](-8,9.5){$E_{\bold{i}}$}
\endpspicture
}
}
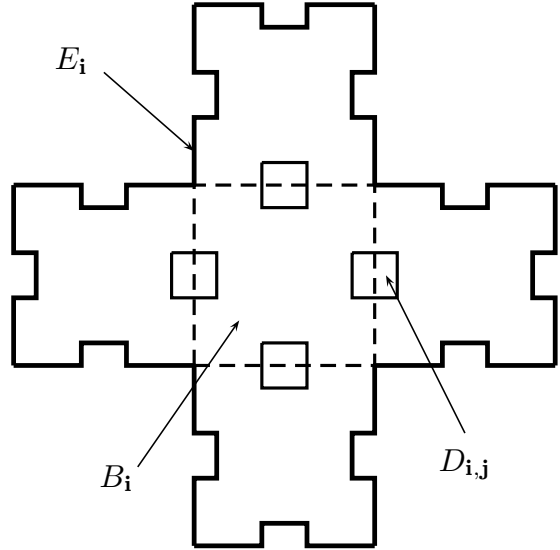
\centerline{figure 3: The region $E_{\bl{i}}$}
\endinsert
\medskip\noindent
The $\sigma$-algebra $\Cal F_\Lambda^{E_\bl{i}}$ is
finer than $\sigma(X_\bl{j}:|\bl{j}-\bl{i}|>1)$, thus it suffices to prove
\eqref{\ineqcond} for $\Phi^\pi[X_\bl{i}=0\;|\;\Cal F_\Lambda^{E_\bl{i}}]$. Clearly $\Cal
F_\Lambda^{E_\bl{i}}$ is atomic and its atoms are of the form $\{\eta\}$, where
$\eta\in\Omega_\Lambda^{E_\bl{i}}$. So let us consider such a
$\eta\in\Omega_\Lambda^{E_\bl{i}}$, then we have that
 $$
\Phi^\pi[X_\bl{i}=0\;|\;\eta]\leq\sum_{\bl{j}\sim\bl{i}}\Phi^{\pi}[K_{\bl{i},\bl{j}}^c\;|\;\eta]+\Phi^{\pi}[R_\bl{i}^c\;|\;\eta].\eqnum
 $$ \eqlabel{\sumcond}\noindent
 For each $\bl{i},\bl{j}\in\Lambda^{(N)}$ such that
$\bl{i}\sim\bl{j}$, let us fix $\eta'\in\Omega_{E_\bl{i}}^{B_{\bl i}}$, $\eta''\in
 \Omega_{E_\bl{i}}^{D_{\bl i, \bl{j}}}$ in order to construct $\eta\eta'\in
 \Omega_{\Lambda}^{B_{\bl i}}$ and $\eta\eta''\in\Omega_{\Lambda}^{D_{\bl i,\bl{j}}}$,
 which are the concatenation of $\eta$ with $\eta'$, respectively with
$\eta''$:
 $$
\eta\eta'(e)=\eta'(e) \hbox{ for } e\in \E(E_\bl{i})\setminus\E(B_\bl{i})\qquad\eta\eta'(e)=\eta(e)\hbox{ for } e\in\E(\Lambda)\setminus \E(E_\bl{i});
 $$
and
 $$
\eta\eta''(e)=\eta''(e) \hbox{ for } e\in \E(E_\bl{i})\setminus\E(D_{\bl{i},
\bl{j}})\qquad \eta\eta''(e)=\eta(e)\hbox{ for } e\in\E(\Lambda)\setminus\E(E_\bl{i}).
 $$
Then, by theorem \procref{\thun}, there exist an integer $N_0>0$ and a real number $C>0$
such that for all $N>N_0$
$$\eqalign{
\Phi^\pi[R_\bl{i}^c\;|\;\eta\eta']&=\Phi^{\pi\cdot
W_\Lambda^{B_\bl{i}}(\eta\eta')}[R_\bl{i}^c]\leq\exp(-C\sqrt{N}),\cr
\Phi^\pi[K_{\bl{i},\bl{j}}^c\;|\;\eta\eta'']&=\Phi^{\pi\cdot
W_\Lambda^{D_{\bl{i},\bl{j}}}(\eta\eta'')}[K_{\bl{i},\bl{j}}^c]\leq\exp(-CN).}$$
Finally, by averaging over all the $\eta'$ and $\eta''$ we get from these
estimates that
$$\eqalign{ \Phi^\pi[X_\bl{i}=0\;|\;\eta]&\leq 4\exp(-C N)+\exp(-C
\sqrt{N})\cr &\leq \exp(-C N^{1/2}),}$$
for $N$ large enough.\qed
 \enddemo
We end this section by proving a useful estimates on the
renormalized process. Let $B(n)$ be a $N$-large box, consider its
$N$-partition and the corresponding $N$-block process. The rescaled
box $B(n)^{(N)}$ will be denoted by $\bl{B}$. For $\delta>0$ we consider the event
$$
Z(n,\delta,N)=\left\{\eqalign{
&\exists!\hbox{ crossing cluster of blocks }\bold{\widetilde{C}}\cr
&\hbox{in }\bold B \hbox{ with }
|\bold{\widetilde{C}}|\geq(1-\delta)|\bold B|
}\right\}.
\eqnum$$\eqlabel{\defZ}
\proclaim{Remark \procnum} The event $Z(n,\delta,N)$
has the following interesting property: the presence of the crossing cluster of blocks
$\bold{\widetilde C}$ induces a set of clusters $\{\widetilde{C}_\bl{i}\hbox{
crossing for }B_\bl{i}:
\bl{i}\in\bold{\widetilde C}\}$ in the original
FK-percolation process. These clusters are connected and form
a crossing cluster $\widetilde{C}$ for $B(n)$.
\endproclaim\proclabel{\GeomZ}
\proclaim{Proposition \procnum}
Let $p>\pgd$ and $q\geq 1$. Then for each $\delta>0$ and $N>0$ large
enough
$$
\limsup_{n\tend\infty}\frac{1}{n}\log\sup_{\Phi\in\Cal R(p,q,B(n))}\Phi\left[Z(n,\delta,N)^c\right]<0.
$$
\endproclaim\proclabel{\EstiRenor}
\demo{Proof}
By theorem 1.1 of \cite{\DeuschPisz}, there exists $p_0\in (0,1)$ such
that for all $p>p_0$,
$$ \limsup_{m\rightarrow \infty} \frac{1}{m}\log\sup_{\und{m}\in\Cal
    H_2(m) }P_{B(\und{m}), \text{site}}^{\,p,\text{ indpt}}\left[
\eqalign{
&\not\exists\hbox{ crossing cluster }\widetilde{C}\hbox{ with}\cr
&|\widetilde{C}| \geq (1-\delta)|B(\und{m})|}
\right]<0.\eqnum
$$\eqlabel{\resDeuschPisz}
\noindent Now choose $N$ such as in proposition \procref{\SM} and such that
$p(N)>p_0$. Then by proposition \procref{\SM} and by
\eqref{\resDeuschPisz} we have that
$$\eqalign{
&\limsup_{n\tend\infty}\frac{1}{n}\log\sup_{\Phi\in\Cal R(p,q,B(n))}\Phi\left[\eqalign{
&\not\exists\hbox{ crossing cluster of blocks }\bold{\widetilde{C}}\cr
&\hbox{in }\bold B \hbox{ with }
|\bold{\widetilde{C}}|\geq(1-\delta)|\bold B|
}\right]\cr
&\cr
\leq&\limsup_{n\rightarrow \infty} \frac{1}{n}\log P_{\bold{B}, \text{site}}^{\,p,\text{ indpt}}\left[
\eqalign{
&\not\exists\hbox{ crossing cluster }\bold{\widetilde{C}}\hbox{ with}\cr
&|\bold{\widetilde{C}}|\geq (1-\delta)|\bold{B}|}\right]<0.\quad\qed
}$$
\enddemo
\subhead{\headnum Proof of the surface order large deviations}\endsubhead
In this section we finally establish theorem
\procref{\prin}. We begin by
stating two lemmas. The first one deals with large deviations
from above. Let $\B(n)$ denote the set of clusters in $B(n)$ intersecting
$\partial B(n)$. Note that if the crossing cluster exists then it is in $\B(n)$.
\proclaim{Lemma \procnum} Let $q\geq 1$ and $p\in[0,1]$.
For $\delta>0$, we have
 $$
 \limsup_{n\tend \infty}\frac{1}{n^2}\log\sup_{\Phi\in \Cal R(p,q,
B(n))}\Phi\left[\sum_{C\in\B(n)}|C|>(\theta+\delta)n^2\right]<0.
 $$
\endproclaim\proclabel{\deviaplus}
 \noindent We omit the proof as it would be an exact repetition of Lemma 5.1 in \cite{\Pisz}.
\par The second lemma is about large deviations from below and is of
surface order, in contrast to lemma \procref{\deviaplus}. In section
 \headref{\connect}, we introduced the event $U(n)=\{\exists! \hbox{
open cluster  }C^* \hbox{ crossing } B(n)\}$.
For $\delta>0$, let us define the event $$ V(n,\delta)=U(n)\cap\{|C^*|>(\theta-\delta)n^2\}.
 $$
 \proclaim{Lemma \procnum} Let $q\geq 1$ and $p>\pgd$. Then for each $\delta>0$,
 $$
 \limsup_{n\rightarrow\infty}\frac{1}{n}\log\sup_{\Phi\in\Cal R(p,q,B(n))}\Phi[V(n,\delta)^c]<0.
\eqnum
 $$
\endproclaim \proclabel{\lemV}
\demo{Proof} For $N>0$, if we set  $\Cal Q(N)=\{x\in B(N), \hbox{dist}(x,\partial B(N))\geq \sqrt{N}\}$ then we have
$$\eqalign{
&\liminf_{n\rightarrow
\infty}\Phi_{B(N)}^f\left[N^{-2}\sum_{C;\text{diam}(C)\geq\sqrt{N}}|C|\right]\cr
\geq&\liminf_{N\tend\infty}N^{-2}\sum_{x\in\Cal
Q(N)}\Phi_{B(N)}^f[\diam(C_x)\geq \sqrt{N}]\cr
\geq&\liminf_{N\tend\infty}N^{-2}\sum_{x\in\Cal
Q(N)}\Phi_{B(x,\sqrt{N})}^f[x\leftrightarrow\partial B(x,\sqrt{N})]\cr
\geq&\liminf_{N\tend\infty}N^{-2}|\Cal Q(N)|\Phi_{B(\sqrt{N})}^f[0\leftrightarrow\partial B(\sqrt{N})]=\theta,
}$$
where the last equality follows from lemma \procref{\thdeux}.
\par\noindent Take $N$ such that
$\Phi_{B(N)}^f[\sum_{C;\text{diam}(C)\geq\sqrt{N}}|C|]\geq(\theta-\delta/4)N^2$,
let $B(n)$ be a $N$-large box and consider its
$N$-partition and the corresponding $N$-block process. The rescaled
box $B(n)^{(N)}$ will be denoted by $\bl{B}$. By proposition \procref{\EstiRenor}, it suffices to give an upper
bound on the probability of the event
 $$
W(n)=Z(n,\delta/8, N)\cap\{|\widetilde{C}|\leq(\theta-\delta)n^2\},
$$
where $N$ is large enough and $Z(n,\delta/8, N)$ is defined in \eqref{\defZ}.
By remark \procref{\GeomZ}, on the event $Z(n, \delta/8, N)$
the crossing cluster $\widetilde{C}$ contains all the $B_\bl{i}$-crossing
clusters $\widetilde{C}_{\bl{i}}$, where $\bl{i}\in\bold{\widetilde C}$ and
$\{B_\bl{i},\ \bl{i}\in\bl{B}\}$ are the partitioning $N$-blocks.
For each $\bl{i}\in\bl{B}$, set $Y_{\bold i}=\sum_{C;\text{diam }C\geq
N^{1/2}}|C|$, where $C$ is a cluster of $B_\bl{i}$. Since
for $\bold i\in
\widetilde{\bold C}$, $Y_{\bold i}=|\widetilde{C}_{\bold i}|,$ we obtain the following lower bound
$$|\widetilde{C}|\geq\sum_{\bold i\in\widetilde{\bold C}}Y_i\geq\sum_{\bold i\in \bold
  B}Y_{\bold i}-\sum_{\bold i\in \bold B\backslash
  \widetilde{\bold C}}|B_{\bold i}|
\geq\sum_{{\bold i\in \bold B \raise 6pt\hbox{$\kern -5.6pt\scriptstyle \circ$}}}Y_{\bold
  i}-(\delta/2)n^2,$$
where  $\bold B\raise 8pt\hbox{$\kern -7.2pt\scriptstyle
\circ$}=\bl{B}\setminus \partial\bl{B}.$\kern
7.2pt Hence on $W(n)$ we have that $\sum_{{\bold i\in \bold B \raise 6pt\hbox{$\kern -5.6pt\scriptstyle \circ$}}}Y_{\bold
  i}\leq (\theta-\delta/2)n^2$.
\par\noindent
Denote by $E(n)$ the event that for each $\bold i\in \bold B\raise
8pt\hbox{$\kern -7.2pt\scriptstyle \circ$}$\kern 7.2pt every edge in
$\partial^{\text{edge}}B_{\bold i}$ is closed. Observing that $\sum_{\bold i\in
  \bold B\raise 6pt\hbox{$\kern -5.6pt\scriptstyle \circ$}}Y_{\bold i}\,$ is an increasing function, we have for each $\Phi\in\Cal R(p,q,B(n)),$
$$
\Phi[W(n)]\leq\Phi^f_{B(n)}\left[\left.\sum_{\bold
i\in \bold B\raise 6pt\hbox{$\kern -5.6pt\scriptstyle \circ$}}Y_{\bold
i}<(\theta-\delta/2)n^2\,\right\vert E(n)\right].
$$
\noindent The variables $(Y_{\bold i},\bold i\in\bold B\raise
8pt\hbox{$\kern -7.2pt\scriptstyle \circ$}$\kern 7.2pt$\!)$ are i.i.d. with respect to the conditional
measure, with an expected value larger than $(\theta-\delta/4)N^2$.
\par\noindent Cram\'er's large deviations theorem yields to
$$
\Phi^f_{B(n)}\left[\left.\frac{1}{n^2}\sum_{\bold
i\in \bold B\raise 6pt\hbox{$\kern -5.6pt\scriptstyle \circ$}}Y_{\bold
i}<\theta-\delta/2\,\right\vert E(n)\right]\leq\exp(-C(\delta, \theta,
N)n^2),
$$
where $C(\delta, \theta,N)$ is a positive constant. This completes the
proof.
\qed\enddemo
\demo{Proof of Theorem \procref{\prin}} First we prove the upper
bound. By lemma \procref{\deviaplus}, we can replace the condition
$n^{-2}|C_m|\in(\theta-\varepsilon, \theta+\varepsilon)$ in the definition of $K(n,
\varepsilon, l)$ by $n^{-2}|C_m|>(\theta-\varepsilon)$ and denote the new but otherwise
unchanged event by $K'(n, \varepsilon, l)$. Set
$$
T(n, \varepsilon, N)=Z(n,
\varepsilon/4, N)\cap\{|\widetilde{C}|>(\theta-\varepsilon)n^2\},
$$
where $Z(n,\varepsilon/4, N)$ is defined by \eqref{\defZ}. Fix $\varepsilon<\theta/2$ and $N$ such as in proposition
\procref{\EstiRenor} and such that  $\sqrt{N}\geq 32/\varepsilon$.
\par\noindent Then by  proposition \procref{\EstiRenor} and by lemma \procref{\lemV},  we have
$$
\limsup_{n\rightarrow\infty}\sup_{\Phi\in \Cal
  R(p,q,B(n))}\frac{1}{n}\log\Phi[T(n, \varepsilon, N)^c]<0.\eqnum$$\eqlabel{\limT}\noindent
Set $n\geq 64 N/\varepsilon$ and $L=2N$, we claim that $T(n, \varepsilon, N)\subset K'(n, \varepsilon, L)$.
This fact, together with \eqref{\limT}, implies the upper
bound. Therefore, to
complete  the upper bound we will proof that the cluster
$\widetilde{C}$ of $T(n, \varepsilon, N)$, is the unique cluster with maximal
volume and that the $L$-intermediate clusters have a negligible
volume. So suppose that  $T(n, \varepsilon, N)$ occurs. As $\varepsilon<\theta/2$
we have that $L^2\leq (\theta-\varepsilon)n^2$, thus the clusters of
diameter less than $L$, have a smaller volume than  $\widetilde{C}$. To control
the size of the clusters different from $\widetilde{C}$ and  of diameter greater than $L$, we define the
following regions:
$$\eqalign{
\forall\,\bold{i}\in \bold{B}\,:\,\quad &G_{\bold i}=\{x\in B_{\bold i}\;|\;\hbox{dist}(x, \partial B_{\bold i})\leq
\sqrt{N}\}\qquad\hbox{and}\qquad \Cal Q_{\bold i}=B_{\bold i}\backslash G_{\bold i},\cr
&\displaystyle G=\bigcup_{\bold i\in\bl{B}}G_{\bold i},
}$$
as shown in figure 4 below:\par
\vbox{
\centerline{
\hbox{
\psset{unit=0.3cm}
\pspicture(-15,-10)(15,11)
\psline[linewidth=1.2pt]{-}(-8.6,-8.6)(8.6,-8.6)
\psline[linewidth=1.2pt]{-}(-8.6,-4)(8.6,-4)
\psline[linewidth=1.2pt]{-}(-8.6,0)(8.6,0)
\psline[linewidth=1.2pt]{-}(-8.6,4)(8.6,4)
\psline[linewidth=1.2pt]{-}(-8.6,8.6)(8.6,8.6)
\psline[linewidth=0.5pt]{<->}(-8.6,9.2)(8.6,9.2)
\uput[90](0,9.2){$n\geq 64N/\varepsilon$}
\psline[linewidth=1.2pt]{-}(-8.6,-8.6)(-8.6,8.6)
\psline[linewidth=1.2pt]{-}(-4,-8.6)(-4,8.6)
\psline[linewidth=1.2pt]{-}(-0,-8.6)(-0,8.6)
\psline[linewidth=1.2pt]{-}(4,-8.6)(4,8.6)
\psline[linewidth=1.2pt]{-}(8.6,-8.6)(8.6,8.6)
\psline[linewidth=0.5pt]{<-}(9,3.5)(9,3.49)
\psline[linewidth=0.5pt]{<-}(9,4.5)(9,4.51)
\psline[linewidth=0.5pt]{-}(9,3.5)(9,4.5)
\uput[0](9,4){$2\sqrt{N}\geq64/\varepsilon$}
\psline[linewidth=0.5pt]{<->}(9,0)(9,-4)
\uput[0](9,-2){$N$}
\pspolygon[fillstyle=vlines](-8.6,-8.6)(8.6,-8.6)(8.6,8.6)(-8.6,8.6)
\psframe[fillstyle=solid](-8.1,-8.1)(-4.5,-4.5)
\psframe[fillstyle=solid](-3.5,-8.1)(-0.5,-4.5)
\psframe[fillstyle=solid](0.5,-8.1)(3.5,-4.5)
\psframe[fillstyle=solid](-8.1,0.5)(-4.5,3.5)
\psframe[fillstyle=solid](-8.1,4.5)(-4.5,8.1)
\psframe[fillstyle=solid](-3.5,4.5)(-0.5,8.1)
\psframe[fillstyle=solid](-8.1,-3.5)(-4.5,-0.5)
\psframe[fillstyle=solid](-3.5,0.5)(-0.5,3.5)
\psframe[fillstyle=solid](-3.5,-3.5)(-0.5,-0.5)
\psframe[fillstyle=solid](4.5,4.5)(8.1,8.1)
\psframe[fillstyle=solid](4.5,0.5)(8.1,3.5)
\psframe[fillstyle=solid](4.5,-8.1)(8.1,-4.5)
\psframe[fillstyle=solid](4.5,-0.5)(8.1,-3.5)
\psframe[fillstyle=solid](0.5,0.5)(3.5,3.5)
\psframe[fillstyle=solid](0.5,4.5)(3.5,8.1)
\psframe[fillstyle=solid](0.5,-3.5)(3.5,-0.5)
\uput[-180](-11,-4){$G_{\bold i}$}
\psline[linewidth=0.5pt]{->}(-11,-4)(-8.2,-3.7)
\uput[-180](-11,-2){$\Cal Q_{\bold i}$}
\psline[linewidth=0.5pt]{->}(-11,-2)(-6.5,-2)
\endpspicture
}
}
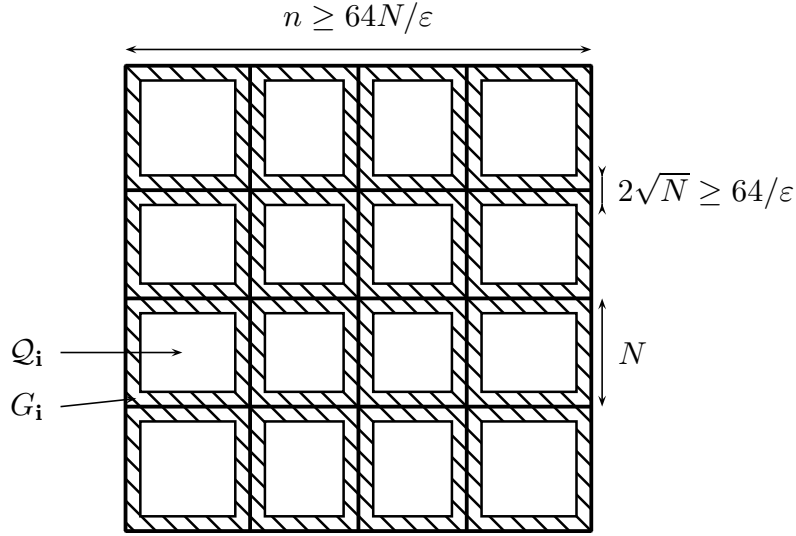
\centerline{figure 4: The regions $G_\bl{i}$ and $\Cal Q_\bl{i}$}
}
\medskip\noindent Then, as $n\geq64N/\varepsilon$, we have
$$
\sum_{\bold{i}\in\partial \bold{B}}|B_{\bold i}|\leq 16 n N\leq\frac{\varepsilon}{4}n^2,
$$
and, as $\sqrt{N}\geq 32/\varepsilon$
$$
|G|\leq 8 \frac{n^2}{\sqrt{N}}\leq\frac{\varepsilon}{4}n^2.
$$
Take a cluster $C$ of diameter greater than $L$ and different from
$\widetilde{C}$. Then $C$ touches at least two
blocks. However, it may not touch the set
$\cup\Cal Q_{\bl{i}}\ $ where $\bl{i}$ runs over $\bold{\widetilde{C}}$; otherwise we would
have that $\diam(C\cap B_\bl{i})\geq\sqrt{N}$ for an occupied block
$B_\bl{i}$, and therefore we would have that $C=\widetilde{C}$.
Hence  all the  clusters of diameter greater than $L$ must
lie in the set  $G\cup(\cup_{\bold i\in \bold{\tilde C}^c}B_{\bold
 i})$. Let us estimate the volume of this set:
$$
|\bigcup_{\bold i\in \widetilde{\bold{C}}^c}B_{\bold i}|\leq\sum_{\bold
 i\in\partial\bold B}|B_{\bold i}|+N^2|\bold{\widetilde
 C}^c|<\frac{\varepsilon}{2}n^2.
$$
Thus
$$
|G\cup(\bigcup_{\bold i\in \widetilde{\bold{C}}^c}B_{\bold i})|\leq\frac{3\varepsilon}{4}n^2.
$$
Since $(3\varepsilon/4)n^2<(\theta-\varepsilon)n^2$, $\widetilde{C}$ is the
unique cluster of maximal volume and the $L$-intermediate class $\J_L$
has a total volume smaller than $(3\varepsilon/4)n^2$. This proves
that $T(n,\varepsilon, L)\subset K'(n,\varepsilon, L)$ and completes
the proof of the upper bound.
\medskip For the lower bound, it suffices to close all the horizontal edges in $B(n)$
 intersecting the vertical line $x=1/2$. This implies that there in no crossing
 cluster in $B(n)$. By \eqref{\NRJ} and FKG inequality, the probability of this event is bounded from below by
 $(1-p)^{n}$.
\qed
\enddemo

\enddocument